\documentclass[12pt]{article}
\usepackage{graphics}
\usepackage{graphicx}
\usepackage{amssymb,amsfonts}
\usepackage{amsmath}
\usepackage{amsthm}

\def\C{\mathbb C}

\newtheorem{theorem}{Theorem}
\newtheorem{result}[theorem]{Result}

\title{Coincidences in numbers of graph vertices corresponding to regular planar hyperbolic mosaics}
\author{L\'aszl\'o N\'emeth\footnote{Institute of Mathematics, University of West Hungary. nemeth.laszlo@emk.nyme.hu}, 
 L\'aszl\'o Szalay\footnote{Institute of Mathematics, University of West Hungary.  szalay.laszlo@emk.nyme.hu}}
 
\date{}

\begin{document}
\maketitle

\def\halmos{\rule{6pt}{6pt}}

\begin{abstract}
The aim of this paper is to determine the elements which are in two pairs of sequences linked to the regular mosaics $\{4,5\}$ and $\{p,q\}$ on the hyperbolic plane. The problem leads to the solution  of diophantine equations of certain types.\\[1mm]
\end{abstract}

\section{Introduction}\label{sec:introduction}

\begin{figure}[!htb]\centering   
\includegraphics{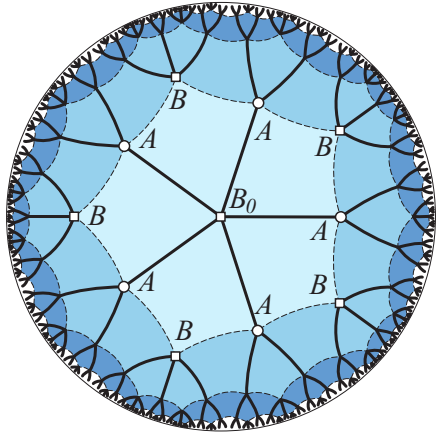} \caption{\emph{Trees of the mosaic $\{p,q\}=\{4,5\}$.}}\label{abra:tree45}
\end{figure}

Consider a regular mosaic on the hyperbolic plane. Such a mosaic is characterized by the Schl\"afli's symbol $\{p,q\}$. It is known that we can define belts of cells around a given vertex of the mosaic (see \cite{nem4}). Let's say that belt ${\cal B}_0$ is the aforesaid fixed vertex itself denoted by $B_0$. The first belt ${\cal B}_1$ consists of the cells which connect  to $B_0$. Assume now that the belts ${\cal B}_{i-1}$ and ${\cal B}_i$ are known ($i\ge1$). Let belt ${\cal B}_{i+1}$ be created by the cells that have common point (not necessarily common vertex) with ${\cal B}_i$, but not with  ${\cal B}_{i-1}$. Figure \ref{abra:tree45} shows the  first three belts in the mosaic corresponding to $\{4,5\}$. One important question is to study the phenomenon of the growing of belts (\cite{hor}, \cite{nem1}, \cite{nem2}), even in higher dimensions, too. 

\begin{figure}[!htb]
\centering   \includegraphics{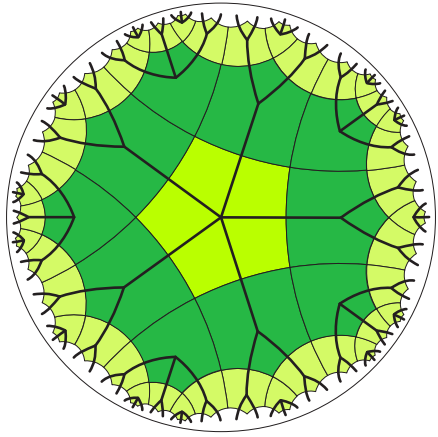} \caption{\emph{Trees of the mosaic $\{5,4\}$, dual of mosaic  $\{4,5\}$.}} \label{abra:tree45dual}
\end{figure}

Take vertex $B_0$ as a main root of a will-be-graph (this is level $0$). In general, let the outer boundary of belt ${\cal B}_i$ be called level $i$. Connect the vertices of level 
$1$ to $B_0$ along the edges between the two levels of the lattice. By this way we have started to build trees. Then use always the maximum number of edges between level $(i-1)$ and  level $i$. All vertices on level $i$ are connected to only one vertex of the previous level, such that no unconnected leaves on level $(i-1)$ are remained. We never connect edges on the same level. The rest vertices on layer $i$ will be roots of new trees. In this way, we obtain infinitely many trees, each of them contains infinitely many vertices. Let $\bar{A}$ denote the set of roots and $\bar{B}$ the set of other vertices. In Figure \ref{abra:tree45} and \ref{abra:tree45dual} the thick edges show the trees from level 0 to level 4. (We remark, that the dual problem is when we establish trees by connecting the centres of the cells of the mosaic.) 

The case $q=3$ provides no any tree since only one edge is not enough to connect the consecutive levels. If $p=3$ the algorithm, apart from $B_0$, does not give roots. 
Therefore we may assume $p\ge4$, $q\ge4$, and since $(p-2)(q-2)=4$ is the Euclidean lattice we also suppose $(p-2)(q-2)>4$.

Let $a_i$ and $b_i$ denote the number of the vertices of $\bar{A}$ and $\bar{B}$ on level $i$, respectively.  
In this paper, we compare the terms $a_i$ (and later $b_i$) of sequences belonging to different Schl\"afli's symbols $\{p,q\}$.

\smallskip

In the following, we recall some properties of the sequences $a_i$ and $b_i$ correspon\-ding to hyperbolic planar lattice $\{p,q\}$ (see \cite{nem4}).
Simple geometric consideration shows $a_1=q$, $b_1=(p-3)q$, further the recursive system
\begin{eqnarray}
a_{n}&=&(q-3)a_{n-1}+(q-2)b_{n-1},\label{mixed}\\
b_{n}&=&\big((q-3)(p-3)-1\big)a_{n-1}+\big((q-2)(p-3)-1\big)b_{n-1}\label{mixed2}
\end{eqnarray}
holds $(n\geq2, p\ge4, q\ge4)$.

It is easy to separate the sequences $\{a_n\}$ and $\{b_n\}$, and it turns out that
\begin{equation}\label{febr}
a_n=\kappa a_{n-1}-a_{n-2}\qquad{\rm and}\qquad b_n=\kappa b_{n-1}-b_{n-2},
\end{equation}
where $\kappa=(p-2)(q-2)-2$ ($\kappa\ge4$). Thus both sequences satisfy the same recurrence relation of order two, and they differ in their initials values. Indeed, to use (\ref{febr}) we need also the terms $a_2$ and $b_2$. Obviously, by (\ref{mixed}) and (\ref{mixed2}), $a_2=(\kappa+1)q$, $b_2=(\kappa(p-3)-1)q$, and $(a_1,a_2)\ne(b_1,b_2)$. Later we also use the term $a_3=(\kappa^2+\kappa-1)q$. Although $a_0$ and $b_0$ have no geometrical meaning, (\ref{febr}) provides the values $a_0=-q$, $b_0=q$, and this sometimes makes the calculations easier.
\smallskip

To achieve the investigations, we introduce the sufficient notations and recall some facts from the theory of linear recurrences. In general, let $r$ and $s$ denote arbitrary complex numbers. The sequence $\{G\}_{n=0}^\infty$ given by the initial values $G_0\in\C$ and $G_1\in\C$, and by the recursive relation
\begin{equation}\label{2olr}
G_n=rG_{n-1}+sG_{n-2}\qquad (n\ge2),
\end{equation} 
is called binary recurrence. For brevity, we often write $G(r,s,G_0,G_1)$ to indicate the parameters of the sequence $\{G\}$.

For any binary recurrence $G(r,s,G_0,G_1)$, the associate sequence of $\{G\}$ is the sequence $H(r,s,H_0,H_1)$ with
\begin{equation}
H_0=2G_1-rG_0\qquad{\rm and}\qquad H_1=rG_1+2sG_0.
\end{equation}

Put $C_G=G_1^2-rG_0G_1-sG_0^2$. It is known that the terms of a binary recurrence $\{G\}$ and its associate sequence $\{H\}$ satisfy the equality
\begin{equation}\label{H2G2}
H_n^2-DG_n^2=4C_G(-s)^n,
\end{equation}
where $D=r^2+4s$.

\section{Preparation and results}

By (\ref{febr}) it follows that the coefficients of the investigated linear recurrences are $r=\kappa$ and $s=-1$. Thus $D=\kappa^2-4$, moreover 
$$
C_a=a_1^2-ra_0a_1-sa_0^2=(\kappa+2)q^2 
$$
and
$$
C_b=b_1^2-rb_0b_1-sb_0^2=((p-3)^2-\kappa(p-3)+1)q^2.
$$
Now we fix a mosaic given by $\{\tilde p,\tilde q\}=\{4,5\}$. Then $\tilde \kappa=4$, $\tilde a_n=4\tilde a_{n-1}-\tilde a_{n-2}$, $\tilde a_1=5$, $\tilde a_2=25$, and
$\tilde b_n=4\tilde b_{n-1}-\tilde b_{n-2}$, $\tilde b_1=5$, $\tilde b_2=15$, moreover $\tilde D=12$.  
The first ten terms of the sequences are given by the following table.

{\renewcommand{\arraystretch}{1.3}
\begin{table}[!ht]
  \centering
\begin{tabular}{|c||c|c|c|c|c|c|c|c|c|c|}
  \hline
 $i$   & 1 & 2  & 3  & 4   & 5   & 6    &  7    & 8     & 9      & 10   \\ \hline \hline
 $\tilde a_i$ & 5 & 25 & 95 & 355 & 1325& 4945 & 18455 & 68875 & 257045 & 959305  \\ \hline
 $\tilde b_i$ & 5 & 15 & 55 & 205 & 765 & 2855 & 10655 & 39765 & 148405 & 553855 \\ \hline
 \end{tabular}
\caption{\emph{Numbers of leaves and roots on level $i$ connected with $\{4,5\}$}\label{tablazat:number_of_vertices}}
\end{table}
}

The associate sequences of $\{\tilde a_n\}$ and $\{\tilde b_n\}$ satisfy
\begin{eqnarray}
\tilde A_n=4\tilde A_{n-1}-\tilde A_{n-2}\qquad{\rm with}\qquad \tilde A_{1}=30,\; \tilde A_{2}=90,\\
\tilde B_n=4\tilde B_{n-1}-\tilde B_{n-2}\qquad{\rm with}\qquad \tilde B_{1}=10,\; \tilde B_{2}=50,
\end{eqnarray}
respectively. Since $C_{\tilde{a}}=150$, $C_{\tilde{b}}=-50$, by (\ref{H2G2}) we obtain the identities
\begin{equation}\label{eqs}
\tilde A_n^2-12\tilde a_n^2=600\qquad{\rm and}\qquad\tilde B_n^2-12\tilde b_n^2=-200.
\end{equation}

\noindent In this paper, we target to solve 

\begin{description}
	\item[I.] the diophantine equation $a_k=\tilde a_\ell$ in $k$ and $\ell$ for certain mosaics $\{p,q\}$ (Section 3); further
	\item[II.] the equations $a_\varepsilon=\tilde a_\ell$ in $\ell$ if $\varepsilon\in\{1,2,3\}$ and one of $p$ and $q$ is fixed (Section 4 and 5).
\end{description}
For the sequence $\{b_n\}$ analogous problems are examined.

The first question leads to simultaneous Pellian equations. The second problem requires different approaches depending on $\varepsilon$ and the sequence $\{a_n\}$ (or $\{b_n\}$).

The observations are contained in the following theorems and Result \ref{th2}. We always assume that 
$$\{p,q\}\ne\{4,4\},\{4,5\}.$$

\begin{theorem}\label{th1}
\begin{enumerate} 
  \item Let $4\le p\le 25$ and $4\le q\le 18$. Then the equation $a_k=\tilde a_\ell$ has only the trivial solution $a_1=\tilde a_1=5$ for $q=5$ and any $p$.

\item If $4\le p,q\le 10$, or $11\le p\le25$ and $4\le q\le 8$, then the equation $b_k=\tilde b_\ell$ possesses only the solutions

\begin{itemize}
   \item $\{p,q\}=\{6,5\}$,  $b_1=\tilde b_2=15$,
   \item $\{p,q\}=\{10,5\}$, $b_2=\tilde b_5=765$,
   \item $\{p,q\}=\{14,5\}$, $b_1=\tilde b_3=55$.
\end{itemize}
\end{enumerate}
\end{theorem}

\begin{result}\label{th2}
\begin{enumerate} 
  \item  If $4\le p\le 1\,600$, then $a_2=\tilde a_\ell$ is satisfied by
    		\begin{itemize}
			      \item  $\{p,q\}=\{26,5\}$,  $a_2=\tilde a_4=335$,
            \item  $\{p,q\}=\{90,29\}$,  $a_2=\tilde a_8=68\,875$,
            \item  $\{p,q\}=\{332,5\}$,  $a_2=\tilde a_6=4\,945$,
		    \end{itemize}
		    	  
  \item  In case of $4\le q\le 10\,000$, $a_3=\tilde a_\ell$ has no non-trivial small solution (i.e.~$p\le10\,000$).
		
  \item Assume $4\le p\le 10\,000$ or $4\le q\le 2\,800$. Then $\{p,q\}=\{10,5\}$,  $b_2=\tilde b_5=765$  
       satisfy the equation $b_2=\tilde b_\ell$.
\end{enumerate}
\end{result}

\begin{theorem}\label{th3}
\begin{enumerate} 
 \item All the solutions to $a_2=\tilde a_\ell$, with $5\le q\le 25$ are given by

\begin{itemize}
   \item $q=5$,  $\ell={2+2t}$ ($t\in\mathbb{N^+}$),
   \item $q=19$,  $\ell={58+90t}$ and $\ell={78+90t}$ ($t\in\mathbb{N}$),
   \item $q=23$,  $\ell={28+88t}$ ($t\in\mathbb{N}$),
   \item $q=25$,  $\ell={32+33t}$ ($t\in\mathbb{N}$).
\end{itemize} 

\item All the solutions to $b_1=\tilde b_\ell$, with $5\le q\le 25$ are given by

\begin{itemize}
   \item $q=9$,   $\ell={5+18t}$ and $\ell={14+18t}$ ($t\in\mathbb{N}$),
   \item $q=11$,  $\ell={3+10t}$ and $\ell={8+10t}$ ($t\in\mathbb{N}$),
   \item $q=15$,  $\ell={2+6t}$ and $\ell={5+6t}$ ($t\in\mathbb{N}$),
	 \item $q=17$,  $\ell={5+18t}$ and $\ell={14+18t}$ ($t\in\mathbb{N}$).
\end{itemize} 
\end{enumerate}
\end{theorem}

\section{Type I:  $a_k=\tilde a_\ell$ and $b_k=\tilde b_\ell$ with certain $p$ and $q$ (Proof of Theorem \ref{th1})}

It is known that the binary recurrence sequences are periodic modulo any positive integer. A simple consideration shows that the terms $\tilde a_n$ are never divisible by 2, 3, 7, 11, 13, 17 (primes up to 25), while $\tilde b_n$ are never a multiple of 2, 7, 13, 19, 23 (primes also up to 25). On the other hand, $q\mid a_n$ and $q\mid b_n$ hold for any $n$. Consequently,  there is no solution to the equation $a_k=\tilde a_\ell$ unless $q=5,19,23,25$. Indeed, by $q\mid a_n$, one needs only to check one period of $\{\tilde a_n\}$ modulo $q$. Similarly, $b_k=\tilde b_\ell$ may possess solution only when $q=5,9,11,15,17,25$.
Unfortunately, we could achive the computations only for $q=5$ regarded to $a_k=\tilde a_\ell$, and for $q=5$ and $q=9$ regarded to $b_k=\tilde b_\ell$ since the time demand of evaluation of the algorithm decribed below seemed to be too much for larger $q$ values. 

Suppose that $p$ and $q$ are given, and consider $a_k=\tilde a_\ell$. Assume that $x=a_k$ satisfies this equation. Then, by (\ref{H2G2})  
\begin{equation}\label{sp1a}
y^2-(\kappa^2-4)x^2=4(\kappa+2)q^2
\end{equation} 
holds for some positive integer $y$. On the other hand, in the virtue of (\ref{eqs}) (the source of (\ref{eqs}) is (\ref{H2G2})), $x=\tilde a_\ell$ is also a zero of the equation
\begin{equation}\label{sp2a}
z^2-12x^2=600
\end{equation} 
for some positive suitable integer $z$.
Clearly,  (\ref{sp1a}) and (\ref{sp2a}) form a system of simultaneous Pellian equations. The {\tt PellianSystem()} procedure, developed in \cite{Sz} and implemented in MAGMA \cite{M} is able to solve such a system if the coefficients are not too large.

If we take $b_k=\tilde b_\ell$, then (\ref{sp1a}) and (\ref{sp2a}) must be replaced by 
\begin{equation}\label{sp1b}
y^2-(\kappa^2-4)x^2=4((p-3)^2-\kappa(p-3)+1)q^2
\end{equation}  
and 
\begin{equation}\label{sp2b}
z^2-12x^2=-200,
\end{equation}
respectively.

We have checked the solutions of the appropriate system of Pellian equations by MAGMA, and the result of the computations is reported in Theorem \ref{th1}.

To illustrate the time demand of the computations, we note that the MAGMA server needed approximately 21 days to show that $b_k=\tilde{b}_\ell$ has no solution in the case $\{p,q\}=\{8,9\}$ (this was the worst case we considered).


\section{Type II: $a_\varepsilon=\tilde a_\ell$, $b_\varepsilon=\tilde b_\ell$, part 1.~(Background behind Result \ref{th2})}\label{hiperell}

This section is devoted to deal with the equations above in the specific cases

\begin{enumerate}
	\item $a_2=\tilde a_\ell$, when parameter $p$ of $\{a_n\}$ is fixed in the range $4\le p\le 1\,600$,
	\item $a_3=\tilde a_\ell$, when parameter $q$ of $\{a_n\}$ satisfies $4\le q\le 10\,000$,
	\item $b_2=\tilde b_\ell$, when $p\in[4;10\,000]$,
	\item $b_2=\tilde b_\ell$, when $q\in[4;2\,800]$.
\end{enumerate}
The common background behind the four problems is that all of them are linked to hyperelliptic diophantine equations of degree four. Observe, that $a_2$ and $b_2$ is a quadratic polynomial in $q$, similarly $a_3$ and $b_2$ has degree two in $p$. 

Consider first $$a_2=\tilde a_\ell$$ with fixed $p$. Then, by the first identity of (\ref{eqs}), $a_2$ satisfies
$$
y^2-12a_2^2=600,
$$  
where $a_2=f(q)=(\kappa+1)q$ is a quadratic polynomial of $q$. Consequently we need to solve the quartic hyperelliptic equation
\begin{equation}\label{he1}
y^2=12f^2(q)+600.
\end{equation}
We use the {\tt IntegralQuarticPoints()} procedure of MAGMA package to handle (\ref{he1}). Note that if the constant term of the polynomial on the right hand side of (\ref{he1}) is not a full square, then the procedure requires a solution (as input) to the equation  to determine all solutions. In this case we scanned the interval $J=[-10\,000;10\,000]$ for $q$  to find a solution. It might occur that there is a solution outside $J$ and not inside $J$, but we found no example to this.

If once we have determined a $q$, then we search back the corresponding subscript $\ell$.

The analogy to the other 3 cases of this section is obvious: in the right hand side of (\ref{he1}) the polyomial $f$ is being replaced by 
$f(p)=(\kappa^2+\kappa-1)q$, $f(q)=(\kappa(p-3)-1)q$ and $f(p)=(\kappa(p-3)-1)q$, respectively.

Solutions we found are listed in Result \ref{th2} (the list might be not full in accordance with the basic interval $J$ which was used for finding a solution).


\section{Type III: $a_\varepsilon=\tilde a_\ell$, $b_\varepsilon=\tilde b_\ell$, part 2.~(Proof of Theorem \ref{th3})}\label{maradekok}

Here we study the title equation in a few cases with small $\varepsilon$, which differ from the previous section. Recall that both of the sequences $\{\tilde a_n\}$ and $\{\tilde b_n\}$ are purely periodic for any positive integer modulus.

Since $a_1=q$ the equation $a_1=\tilde a_\ell$ has, trivially, infinitely many solutions.
\smallskip

The next problem is $a_2=\tilde a_\ell$ with fixed $q$. (The case with fixed $p$ has already been studied in Section \ref{hiperell}.) Recall that $a_2=(\kappa+1)q$, more precisely $$a_2=q(q-2)(p-2)-q$$ 
is linear in $p$. Therefore we need to determine the common terms of an arithmetic progression and the sequence $\{\tilde a_n\}$. The situation does not change if we consider $b_1=\tilde b_\ell$ with either fixed $p$ or fixed $q$. Indeed, $b_1=(p-3)q$ is linear both in $p$ and $q$.

Obviously, $a_2\equiv -q\;(\bmod \,\,q(q-2))$. Consequently, the equation $a_2=\tilde a_\ell$ is soluble if and only if we find at least one element in the sequence $\{\tilde a_n\}$, which is congruent $-q$ modulo $q(q-2)$. Because of the periodicity, one must check only one period of $\{\tilde a_n\}$ modulo $q(q-2)$. 

Assume first that $q=5$. Then for the modulus $q(q-2)=15$ we have $\tilde a_{2+2t}\equiv -5$ (the cycle's length is 2, and $t\in\mathbb{N}$). Hence $a_2=\tilde a_{2+2t}$, further
$$
p=\frac{\tilde a_{2+2t}+q}{q(q-2)}+2.
$$
The first six $t$ values yield the following solutions. (If $t=0$ then the two sequences $\{a_n\}$ and $\{\tilde a_n\}$ coincide.)
{\renewcommand{\arraystretch}{1.3}
\begin{table}[!ht]
  \centering
\begin{tabular}{|l||c|c|c|c|c|c|}
  \hline
 $t$   & 0 & 1  & 2  & 3   & 4  & 5    \\ \hline \hline
 $a_2=\tilde a_{2+2t}$ & 25 & 355 & 4945 & 68875 & 959305 & 13361395  \\ \hline
 $p$ & 4 & 26 & 332 & 4594 & 63956 & 890762\\ \hline
 \end{tabular}
\caption{\emph{First few solutions to $a_2=\tilde a_\ell$ when $q=5$}\label{tablazat:uto}}
\end{table}
}

If $q>5$ the first non-trivial solution is occurred when $q=19$. Here the length of the cycle is 90, and $q(q-2)\mid\tilde a_{58}+19$, $q(q-2)\mid \tilde a_{78}+19$. That is $a_2=\tilde{a}_{58+90t}$ and $a_2=\tilde{a}_{78+90t}$ ($t\in\mathbb{N}$) provide all solutions for suitable values $p$. For instance, $t=0$ gives 
$$p=8\,437\,940\,669\,128\,098\,583\,408\,551\,589\,590$$ 
and 
$$p=2\,318\,394\,927\,973\,629\,460\,854\,893\,981\,169\,574\,319\,067\,870,$$ 
respectively.

The treatment is similar for $b_1=\tilde b_\ell$. If $q=5$, then solution always exists since $b_1=(p-3)q$, $5\mid \tilde b_\ell$, therefore $p=\tilde b_\ell/5+3$. ($\tilde b_2$ and $\tilde b_3$ give back solutions have already been appeared in Theorem \ref{th1}.) Now $b_1\equiv 0\;(\bmod \,q)$, and fixing $q\ge6$ the first solution appears for $q=9$, when the cycle length is 18 (modulo $q$), and we have $b_1=\tilde{b}_{5+18t}$ and $b_1=\tilde{b}_{14+18t}$ ($t\in\mathbb{N}$). These results can  be directly converted the results corresponding to $p$, therefore we omit the appropriate analysis. 

The results we obtained are summarized in Theorem \ref{th3}.

\smallskip

Finally, we examine the equation $a_3=\tilde a_\ell$ with fixed $q$, further $b_3=\tilde b_\ell$ when exactly one of $p$ and $q$ is given. In each case we have a polynomial of degree three, let say $\phi(x)$, and we look for the common values of the polynomial and a binary recurrence. By (\ref{H2G2}), the problem leads to the hyperelliptic equation 
$$
y^2=12\phi^2(x)+c
$$
of degree 6, where the constant $c$ is either $600$ or $-200$. Since the leading coefficient on the right hand side is not a square, there is no genearal algorithm to solve.
For example, $p=5$ provides now 
$$
y^2=12q^2(9q^2-45q+55)^2+600.
$$
After dividing by 4, we have 
$$
y_1^2=243q^6-2430q^5+9045q^4-14850q^3+9075q^2+150,
$$
and the techique of the solution is not known.

\medskip

{\bf Acknowledgements.} The authors thank P.~Olajos for his valuable help in using MAGMA package.

\end{document}